\documentclass[12pt,draft]{article}
\usepackage[cp1251]{inputenc}
\usepackage[english,russian]{babel}
\usepackage{amsmath,amscd,amssymb,cite}
\begin{document}

\newtheorem{theorem}{Теорема}
\newtheorem{lemma}{Лемма}
\newtheorem{proposition}{Предложение}
\newtheorem{axiom}{Аксиома}
\newtheorem{example}{Пример}
\newtheorem{remark}{Замечание}
\newtheorem{state}{Утверждение}

\def \Pil{\mathop{\Pi}\limits}
\def \vl{\mathop{\vee}\limits}
\def \liml{\mathop{\lim}\limits}
\def \supl{\mathop{\sup}\limits}
\def \infl{\mathop{\inf}\limits}
\def \cupl{\mathop{\cup}\limits}
\def \intl{\mathop{\int}\limits}
\def \suml{\mathop{\sum}\limits}
\def \mil{\mathop{\min}\limits}
\def \mal{\mathop{\max}\limits}
\def \ll{\mathop{\&}\limits}
\def \ulil{\mathop{\underline{\lim}}\limits}
\def\huh{\hbox{\vrule width 2pt height 8pt depth 2pt}}
\def\op{\begin{theorem} \rm}
\def\Bbb#1{{\bf #1}}
\def\ddd#1#2{buildrel \#1 \#2}
\def\eqnum#1{\eqno (#1)}
\def\fnote#1{\footnote}
\def\blacksquare{\hbox{\vrule width 4pt height 4pt depth 0pt}}
\def\square{\hbox{\vrule\vbox{\hrule\phantom{o}\hrule}\vrule}}

\newcommand{\I}{{\rm I}\hspace{-0.2em}{\rm I}}
\newcommand{\R}{{\mathbb R}}
\newcommand{\kb}[1]{{\mathbb{#1}}}
\newcommand{\A}{\protect\text{\boldmath$A$}}
\newcommand{\ab}{\protect\text{\boldmath$a$}}
\newcommand{\bb}{\protect\text{\boldmath$b$}}
\newcommand{\cb}{\protect\text{\boldmath$c$}}
\newcommand{\db}{\protect\text{\boldmath$d$}}
\newcommand{\ub}{\protect\text{\boldmath$u$}}
\newcommand{\dual}{{\rm dual}\;}
\newcommand{\ov}{\overline}
\newcommand{\un}{\underline}
\newcommand{\ost}{\odot\hspace{-3.6mm}\star \;}
\newcommand {\rx}{{\rm x}}
\newcommand {\rw}{{\rm w}}
\newcommand {\ru}{{\rm u}}
\newcommand {\ry}{{\rm y}}
\newcommand {\rp}{{\rm p}}
\newcommand {\rz}{{\rm z}}

\newcommand{\mbf}[1]{\protect\text{\boldmath$#1$}}

\begin{center}
БЕСКВАНТОРНОЕ ОПИСАНИЕ  МНОЖЕСТВА РЕШЕНИЙ ОБОБЩЕННОЙ \\
ИНТЕРВАЛЬНО-КВАНТОРНОЙ СИСТЕМЫ ЛИНЕЙНЫХ УРАВНЕНИЙ \\ 
А.В. Лакеев \\ 
Институт динамики систем и теории управления имени В.М. Матросова СО РАН,\\ г. Иркутск, Россия  \\
e-mail: lakeyev@icc.ru \\
И.А. Шарая \\ 
Институт вычислительных технологий СО РАН, г. Новосибирск, Россия \\
e-mail: sharaya@ict.nsc.ru
\end{center}

\begin{center}
{\bf  Введение}
\end{center}

Многие задачи, возникающие при исследовании реальных систем, содержат не\-точ\-ные данные. Одним из способов описания таких данных является задание интервалов, в которых они могут изменяться. Такой подход привел к созданию в 60-х годах прошлого столетия интервального анализа, интенсивное развитие которого продолжается и в настоящее время.

В данной работе мы будем рассматривать следующую интервальную систему линейных алгебраических  уравнений (ИСЛАУ) (и некоторые ее обобщения):
$$
\A x =\bb,
\eqno (1)
$$
где $\A \in \kb{IR}^{m\times n}$ — вещественная интервальная $m\times n$-матрица, $\bb \in \kb{IR}^m$ -- $m$-мерный 
интервальный вектор, $x\in \kb{R}^n$. При этом будем в основном придерживаться терминологии и обозначений из \cite{Ker, sarb}.

Множество решений такой системы  может быть определено различными способами в зависимости от того, какими кванторами связываются коэффициенты матрицы и правой части \cite{vat, sarb, irsar, irsar1}. 

В качестве частных  случаев множеств решений ИСЛАУ
рассмотрим следующие.

Исторически первым и до сих пор самым изучаемым является так называемое 
объединенное множество решений ИСЛАУ ({\em united solution set}) 
(W. Oettli,  W. Prager (1964) \cite{oetpr, Oet}, 
H. Beeck (1972) \cite{Beeck}, J. Rohn (1984) \cite{Rohn1}, С.П. Шарый (1990) \cite{sar1}) вида
$$\varXi_{uni}(\A, \bb) = \{ x \in \R^n \; |
  \; \exists A \in \A \;\; \exists b \in \bb \;\; Ax = b \}.
$$
В этом определении все коэффициенты матрицы и правой части связаны кванторами существования $\exists$.
                                                                                       
Затем из практических задач появились и стали изучаться и некоторые другие множества решений
системы:

-  допусковое множество решений ИСЛАУ ({\em tolerable solution set}) 
(E. Nuding, 
\linebreak
J. Wilhelm (1972) \cite{nud}, J. Rohn (1985) \cite{Rohn3}, A. Neumaier (1986) \cite{Neum1}, В.В. Шайдуров, С.П. Шарый (1988) \cite{sar2}, С.П. Шарый (1991) \cite{sar3})
$$
 \varXi_{tol}(\A, \bb) = \{ x \in \R^n \; | \;
 \forall A \in \A \;\;  \exists b \in \bb \;\; Ax = b \},
$$
которое возникло при решении задачи о допусках  \cite{sar3};

- управляемое множество решений ИСЛАУ ({\em controllable solution set}) 
(Н.А. Хлебалин, 
Ю.И. Шокин (1991) \cite{xsok}, Лакеев А.В., Носков С.И. (1992) \cite{ln},  С.П. Шарый (1992) \cite{sar3a, sar3b})
$$
 \varXi_{ctr}(\A, \bb) = \{ x \in \R^n \; | \;
 \forall b \in \bb \; \; \exists A \in \A \; \; Ax = b \},
$$
которое возникло при решении задачи автоматического регулирования
в интервальной постановке \cite{xsok} и  задачи управляемости для системы типа вход-выход \cite{sar3a, sar3b}.

Все эти множества являются частными случаями AE-решений, введенных 
\linebreak
С.П. Шарым (1995)  \cite{sar4, sar5}, в которых все кванторы всеобщности предшествуют кванторам существования. Определим AE-решения (AE-solution set) следуя 
\cite[Глава 5.2]{sarb}.

Пусть для уравнения (1) вместе с интервальной матрицей $\A$ и интервальным вектором $\bb$ заданы $m\times n$-матрица $\alpha$ и  $m$-мерный  вектор $\beta$, состоящие из кванторов. 

Определим интервальные матрицы   
${\A}^{\forall} = ({\ab}_{ij}^{\forall})$, 
${\A}^{\exists} = ({\ab}_{ij}^{\exists})$ и интервальные векторы $\bb^{\forall} =(\bb_i^{\forall})$,
$\bb^{\exists} =(\bb_i^{\exists})$  
тех же размеров, что $\A$ и $\bb$ соответственно, следующим образом:
$$
{\ab}_{ij}^{\forall}
=
\left\{
\begin{array}{rl}
{\ab}_{ij}, & \mbox{если} \ \alpha_{ij} = \forall, \\
0, &  \mbox{если} \ \alpha_{ij} = \exists,  
\end{array} \right.
\hspace{22mm}
{\ab}_{ij}^{\exists}
=
\left\{
\begin{array}{rl}
0, & \mbox{если} \ \alpha_{ij} = \forall, \\
{\ab}_{ij}, &  \mbox{если} \ \alpha_{ij} = \exists,  
\end{array} \right.
$$
$$
{\bb}_{i}^{\forall}
=
\left\{
\begin{array}{rl}
{\bb}_{i}, & \mbox{если} \ \beta_{i} = \forall, \\
0, &  \mbox{если} \ \beta_{i} = \exists,  
\end{array} \right.
\hspace{22mm}
{\bb}_{i}^{\exists}
=
\left\{
\begin{array}{rl}
0, & \mbox{если} \ \beta_{i} = \forall, \\
{\bb}_{i}, &  \mbox{если} \ \beta_{i} = \exists.  
\end{array} \right.
$$
Тогда очевидно, что $\A = {\A}^{\forall} + {\A}^{\exists}$, $\bb = {\bb}^{\forall} + {\bb}^{\exists}$ и выполняются соотношения  дизъюнктности ${\ab}_{ij}^{\forall} {\ab}_{ij}^{\exists} =0$, ${\bb}_{i}^{\forall} {\bb}_{i}^{\exists} =0$. 

{\bf Определение 1} \cite[Определение 5.2.2]{sarb}. Пусть для системы (1) заданы кванторные матрица $\alpha$ и 
вектор $\beta$. 
Множеством AE-решений типа $\alpha\beta$ для (1) называется множество
$$
\varXi_{\alpha\beta}(\A, \bb) = \{ x \in \R^n \; |
\; \forall \hat{A} \in {\A}^{\forall}  \;\; \forall \hat{b} \in {\bb}^{\forall} \;\;
\exists \breve{A} \in {\A}^{\exists} \;\; \exists \breve{b} \in {\bb}^{\exists}
\;\; (\hat{A} + \breve{A})x=\hat{b}+\breve{b}\}.
$$ 

Для этого множества в \cite{sar4} получено следующее бескванторное описание.

{\bf Теорема 1} (характеризация Шарого множеств AE-решений \cite[Теорема 5.2.2]{sarb}). 
Точка $x \in \R^n$ принадлежит множеству AE-решений $\varXi_{\alpha\beta}(\A, \bb)$ тогда и
только тогда, когда
$$
{\A}^{\forall}\cdot x - {\bb}^{\forall} \subseteq {\bb}^{\exists} - {\A}^{\exists}\cdot x ,
\eqno (2)
$$
где « · » — интервальное матричное умножение.

Следующее бескванторное описание ввиде разрешимости системы линейных неравенств с модулями получено в \cite{R1, R2}.

{\bf Теорема 2} (характеризация Рона множеств AE-решений \cite[Теорема 5.2.4]{sarb}). Точка $x \in \R^n$ принадлежит множеству 
AE-решений $\varXi_{\alpha\beta}(\A, \bb)$ тогда и
и только тогда, когда
$$
|({\rm mid}(\A)\cdot x - {\rm mid}(\bb)|
\leq 
({\rm rad}(\A^{\exists}) - {\rm rad}(\A^{\forall}))\cdot |x| + 
{\rm rad}(\bb^{\exists}) - {\rm rad}(\bb^{\forall}),
\eqno (3)
$$
где ${\rm mid}(\cdot)$ и ${\rm rad}(\cdot)$ -- обозначения для середины и радиуса, соответственно (интервальных матриц и векторов \cite{sarb}), а $|x|$ понимается покоординатно.

В связи с понятием множества AE-решений И. Рон \cite{R1} предложил рассматривать несколько более общее уравнение, чем (1). 

Будем считать, что заданы две интервальные матрицы 
$\A^{\prime}, \, \A^{\prime\prime}$ и два интервальных вектора $\bb^{\prime}, \, \bb^{\prime\prime}$. Рассмотрим следующее уравнение: 
$$
(\A^{\prime} + \A^{\prime\prime}) x = \bb^{\prime} + \bb^{\prime\prime}.
\eqno (4)
$$ 
При этом матрицы и векторы -- произвольные (не связаные соотношениями дизъюнктности). 

{\bf Определение 2} \cite{R1}.
Множеством решений уравнения (4) называется множество
$$
\begin{array}{r}
\varXi_{\forall\exists} (\A^{\prime}, \A^{\prime\prime}, \bb^{\prime}, \, \bb^{\prime\prime}) =
\{x\in \R^n \; | \; \forall A^{\prime} \in
\A^{\prime} \;\; \forall b^{\prime} \in \bb^{\prime} \\
$$ 
$$
\exists A^{\prime\prime} \in \A^{\prime\prime} \;\; \exists b^{\prime\prime} \in \bb^{\prime\prime}
\;\; (A^{\prime}+A^{\prime\prime})x=b^{\prime}+b^{\prime\prime}\}.
\end{array}
\eqno (5)
$$ 

И. Рон \cite{R1} заметил, что для множества 
$\varXi_{\forall\exists} (\A^{\prime}, \A^{\prime\prime}, \bb^{\prime}, \, \bb^{\prime\prime})$ также верна характеризация Шарого, то есть формула $(2^{\prime})$, аналогичная (2):
$$
x \in \varXi_{\forall\exists} (\A^{\prime}, \A^{\prime\prime}, \bb^{\prime}, \, \bb^{\prime\prime}) \Leftrightarrow
\A_1^\prime \cdot x - \bb_1^\prime
\subseteq 
\bb_1^{\prime\prime} - \A_1^{\prime\prime} \cdot x ,
\eqno (2^{\prime})
$$
и характеризация в виде разрешимости системы линейных неравенств с модулями,  то есть формула $(3^{\prime})$, аналогичная (3):
$$
x \in \varXi_{\forall\exists} (\A^{\prime}, \A^{\prime\prime}, \bb^{\prime}, \, \bb^{\prime\prime}) \Leftrightarrow
|({\rm mid}(\A^{\prime}) +{\rm mid}(\A^{\prime\prime}))\cdot x - ({\rm mid}(\bb^{\prime}) +{\rm mid}(\bb^{\prime\prime}))|
\leq
$$
$$ 
\leq 
({\rm rad}(\A^{\prime\prime}) - {\rm rad}(\A^{\prime}))\cdot |x| + 
{\rm rad}(\bb^{\prime\prime}) - {\rm rad}(\bb^{\prime}).
\eqno (3^{\prime})
$$

{\bf Замечание 1.} Очевидно, что если для уравнения (4) выполняются условия 
дизъюнктности ${\ab}_{ij}^{\prime} {\ab}_{ij}^{\prime\prime} =0$, ${\bb}_{i}^{\prime} {\bb}_{i}^{\prime\prime} =0$, то, 
полагая $\A =\A^{\prime} +\A^{\prime\prime}$, $\bb =\bb^{\prime} +\bb^{\prime\prime}$ практически однозначно 
(с точностью до нулевых элементов $\A$ и $\bb$, которые не влияют на множество решений),
находим матрицу $\alpha$ и вектор $\beta$ такие, что 
$\varXi_{\alpha\beta}(\A, \bb) = 
\varXi_{\forall\exists} (\A^{\prime}, \A^{\prime\prime}, \bb^{\prime}, \, \bb^{\prime\prime})$.

Опираясь на формулы (3) и $(3^{\prime})$, нетрудно показать, что хотя уравнение (4) и выглядит более общим, чем уравнение (1), но класс его решений вида (5) совпадает с классом AE-решений для уравнения (1). 
В действительности это является следствием того, что в виде неравенства (3) можно представить множество решений  любой  системы линейных неравенств с модулями вида 
$$
|Cx-c|\leq D|x|+d,
\eqno (lm)
$$ 
где $C,\;D$ -- $m\times n$-матрицы, $c,\;d$ -- $m$-мерные векторы, $x\in \kb{R}^n$.
Более точно верно следующее утверждение, которое понадобится нам и в дальнейшем.

{\bf Предложение 1.} Для любых  матриц 
$C,\;D$ и  векторов $c,\;d$ существуют интервальная 
матрица $\tilde{\A}$, интервальный вектор $\tilde{\bb}$ и кванторные матрица $\alpha$ и вектор $\; \beta$ такие, что выполняется равенство
$$\varXi_{\alpha\beta}(\tilde{\A}, \tilde{\bb}) =
\{x\in \R^n \; | \; |Cx-c|\leq D|x|+d \; \}. 
$$

{\bf Доказательство.} Представляя интервальную матрицу $\tilde{\A}$ и интервальный вектор $\tilde{\bb}$ в
центрально-симметричной форме, полагаем 
$$
\tilde{\A} = [C-|D|, \; C+|D|], \;\;\; 
\hspace{10mm}
\tilde{\bb} = [c-|d|, \; c+|d|],
$$
$$
\alpha_{ij} =
\left\{
\begin{array}{ll}
\exists, & \mbox{если} \ d_{i,j} \geq 0 , \\
\forall, &  \mbox{если} \ d_{i,j} < 0 ,
\end{array} \right. \;
\hspace{10mm}
\beta_{i} =
\left\{
\begin{array}{ll}
\exists, & \mbox{если} \  d_{i} \geq 0 , \\
\forall, &  \mbox{если} \  d_{i} < 0 . \\
\end{array} \right.
$$
Для так определенных матриц и векторов, очевидно, выполняются следующие равенства:
$$
{\rm mid}(\tilde{\A}) = C, \;\;\;
{\rm rad}(\tilde{\A}^{\exists}) = D^+ , \;\;\;
{\rm rad}(\tilde{\A}^{\forall}) = D^- ,
$$
$$
{\rm mid}(\tilde{\bb}) = c, \;\;\;
{\rm rad}(\tilde{\bb}^{\exists}) = d^+ ,\;\;\;
{\rm rad}(\tilde{\bb}^{\forall}) = d^- ,
$$
где $a^+$ и $a^-$ -- положительная часть и отрицательная часть числа $a$,  соответственно, а для матриц и векторов 
положительная и отрицательная части
определяются поэлементно \cite[стр. 52]{sarb}. 
\newline
\indent
Подставляя эти равенства в формулу (3) и учитывая, что $a^+ -a^- =a$, получаем неравенство $(lm)$ и, следовательно, равенство 
$\{x\in \R^n \; | \; |Cx-c|\leq D|x|+d \; \} =$ 
\linebreak 
$=\varXi_{\alpha\beta}(\tilde{\A}, \tilde{\bb})$. 
Предложение доказано.

{\bf Следствие 1.} Для любых интервальных матриц 
$\A^{\prime}, \, \A^{\prime\prime}$ и интервальных векторов $\bb^{\prime}, \, \bb^{\prime\prime}$ существуют интервальная 
матрица $\tilde{\A}$, интервальный вектор $\tilde{\bb}$ и кванторные матрица $\alpha$ и вектор $\; \beta$ такие, что выполняется равенство
$$\varXi_{\alpha\beta}(\tilde{\A}, \tilde{\bb}) = 
\varXi_{\forall\exists} (\A^{\prime}, \A^{\prime\prime}, \bb^{\prime}, \, \bb^{\prime\prime}).
$$
Для доказательства достаточно положить $C= {\rm mid}(\A^{\prime}) +{\rm mid}(\A^{\prime\prime})$, 
$D= {\rm rad}(\A^{\prime\prime}) -$ $- {\rm rad}(\A^{\prime})$,  
$c= {\rm mid}(\bb^{\prime}) +{\rm mid}(\bb^{\prime\prime})$ и 
$d= {\rm rad}(\bb^{\prime\prime}) - {\rm rad}(\bb^{\prime})$ в предложении 1. 

Особо отметим работу А.А. Ватолина (1984) \cite{vat}, посвященную задачам линейного программирования с
интервальными коэффициентами. В этой работе 
было введено понятие множества решений
интервальной системы с произвольной кванторной приставкой, однако с дополнительным ограничением неотрицательности 
на переменные $x\in \kb{R}^n$. Такое ограничение существенно упрощает задачу описания этих множеств. В частности, 
именно за счет неотрицательности $x$  А.А. Ватолин получил описание множеств решений в виде разрешимости некоторой 
системы линейных неравенств, что совершенно не характерно для интервальных систем. 

В этой работе также имеется и интерпретация таких множеств решений,  
как решения некоторых игр или многошаговых процессов принятия решений
в условиях интервальной неопределенности, т.е. как решения минимаксных задач исследования операций.
Более подробное описание этой интерпритации имеется в монографии 
С.П. Шарого \cite[стр. 204]{sarb}.
Отметим, что отсутствие хорошего бескванторного описания этих множеств, по-видимому, является одной из причин того, что 
они пока не получили широкого применения. 

В работе  И.А. Шарой (2014) \cite{irsar, irsar1} было введено понятие интервально-кванторных {\bfуравнений и неравенств}, 
которое включает в себя все предыдущие. Кроме того, в ней  
получено \cite[Теорема 2]{irsar} полное описание множества решений интервально-кван\-тор\-ной линейной системы 
{\bf неравенств}. 
Поэтому в данной работе мы будем рассматривать только интервально-кван\-тор\-ные линейные системы {\bf уравнений}.

Задача бескванторного описания различных  множеств решений системы линейных интервальных уравнений имеет давнюю историю.
Первый результат был получен Оеттли и Прагером в 1964 году (W. Oettli,  W. Prager (1964) \cite{oetpr, Oet}). 

Обычно стараются перейти к описанию в вещественной арифметике 
\cite[гл. 2]{FNR},  \cite[c. 93–95]{Ere},  \cite{oetpr, Oet, ger, vat, Rohn3, ln, lak1, lak2, lak3, sar5a, sar5b},
поскольку она привычна, обладает хорошими свойствами и развитыми численными
методами. Ряд бескванторных описаний получен в интервальных арифметиках для различных подклассов интервально-кванторных систем линейных уравнений \cite{sarb, sar4, Beeck1, Neum1}, 
и несмотря на плохие свойства этих арифметик (отсутствие дистрибутивности и т.п.), найденные описания оказались полезны. Так, описание множеств AE-решений интервальных систем линейных
уравнений позволило построить теорию этих множеств и интервальные методы их оценивания
(например, интервальный метод Гаусса—Зейделя и формальный алгебраический подход) \cite{sarb, sar6}.

Целью данной работы является  описание множества решений интервально-кван\-тор\-ной системы линейных  уравнений  
и некоторого его обобщения как в интервальной арифметике, так и 
в виде разрешимости систем  линейных неравенств с модулями, обобщая, таким
образом, теоремы 1 и 2. 
Это позволяет применять уже созданные методы (в том числе и численные) для их исследования.

\begin{center}
{\bf 1.  Интервально-кванторные системы   \\ линейные уравнений в $\kb{IR}^n$ и их обобщение }
\end{center}

Интервалом в классической интервальной арифметике $\kb{IR}$ называют непустое ограниченное связное замкнутое 
подмножество числовой оси. Согласно стандарту на обозначения \cite{Ker} интервальные объекты, в отличие от 
точечных (неинтервальных), будем выделять жирным шрифтом.

Следуя работе  И.А. Шарой  \cite{irsar, irsar1}, введем понятие интервально-кванторной системы линейных уравнений следующим образом.

Рассмотрим систему линейных уравнений вида
$$
Ax=b, \;\; A\in \kb{R}^{m\times n}, \;\; x\in \R^n , \;\; b\in \R^m, 
$$
где $x$ — вектор неизвестных, 
а всякий параметр $u\in \R$ 
(элемент матрицы $A$ или правой части $b$) может принимать значение в пределах
заданного одноименного интервала $\ub$ из $\kb{IR}$. С каждым параметром $u$ свяжем квантор всеобщности либо существования и соответствующую элементарную кванторную приставку  ($\forall u\in \ub$)
либо ($\exists u \in \ub$). Такую интервальную неопределенность параметров можно задать интервальной матрицей 
$\A \in \kb{IR}^{m\times n}$, матрицей кванторов $\cal A$ тех же размеров, что и $\A$, интервальным
вектором $\bb \in \kb{IR}^m$ и вектором кванторов $\beta$ длины $m$. Запишем все элементарные кванторные
приставки в произвольном порядке и обозначим полученную приставку длины $m(n + 1)$ как $Q(\A,\bb, \cal A, \beta )$.

{\bf Определение 3} (И.А. Шарая  \cite{irsar, irsar1}). Интервально-кванторной системой линейных 
уравнений
будем называть формулу
$Q(\A,\bb, {\cal A}, \beta ) (Ax = b)$, а ее решением — всякий вектор $x \in \R^n$, для которого 
формула принимает значение “истина”.

В этом определении присутствует кванторная приставка $Q(\A,\bb, \cal A, \beta )$, которая плохо структурирована. Поэтому естественно возникает желание привести ее к некоторому виду, который был бы  
более конкретным, но так, чтобы  при этом не потерять общность.

Оказалось, что это можно сделать так же, как и для  AE-решений, с помощью некоторых дизъюнктных разбиений матрицы $\A$ 
и вектора $\bb$.

Вначале представим кванторную  приставку $Q(\A,\bb, \cal A, \beta )$ более детально. 
Обозначим через $\left( A \; b\right) \;\;\; m\times (n + 1)$-матрицу, образованную приписыванием справа к матрице $A$ столбца $b$. 
Пусть $\mu = m(n + 1)$ и кортеж 
$$T=\left\langle a_{11},\ldots ,a_{1n},b_1, a_{21},\ldots ,a_{2n},b_2,\ldots ,a_{m1},\ldots ,a_{mn},b_m\right\rangle$$
составлен из элементов матрицы $\left( A \; b\right)$, выписанных в данном порядке. 

Пусть $Q(\A,\bb, {\cal A}, \beta )=\Omega_{\mu} \ldots \Omega_2\Omega_1$, где 
$\Omega_i \in \left\{ (\exists u_i \in \ub_i) , (\forall u_i \in \ub_i)\right\}$ для всех 
\linebreak
$i=\overline{1,\mu}$, и 
кортеж 
$\left\langle u_1, \ldots , u_{\mu}\right\rangle$ 
является перестановкой кортежа $T$.

Разобъем  $Q(\A,\bb, {\cal A}, \beta )$ на $\kappa$ блоков (AE-блоки), т.е. представим в виде
\linebreak 
$Q(\A,\bb, {\cal A}, \beta )=B_{\kappa} \ldots  B_1$, где $B_s = \Omega_{i_s} \ldots \Omega_{i_{s-1}+1}$ для $s=\overline{1,\kappa}$ 
($i_0 =0, \; i_{\kappa}=\mu$) так, чтобы выполнялись следующие условия:

1) для $s=\overline{2,\kappa -1}$,  $\Omega_{i_s}$ -- квантор всеобщность, $\Omega_{i_{s-1}+1}$ -- квантор существования и 
внутри $B_s$ происходит единственая смена смысла кванторов;

2) $B_1$ либо  удовлетворяет условию 1), либо   состоит только из кванторов всеобщности;

3) $B_{\kappa}$ либо  удовлетворяет условию 1), либо  состоит только из кванторов существовония.

Очевидно, что условиями 1)--3) разбиение $Q(\A,\bb, {\cal A}, \beta )$ на AE-блоки определяется однозначно.

{\bf Определение 4}. Для заданных матрицы $\A$, вектора $\bb$ и кванторной приставки $Q(\A,\bb, \cal A, \beta )$ определим
кортежи  
$\A_\forall =\left\langle \A_1^\prime , \ldots ,\A_{\kappa}^\prime \right\rangle$, 
$\A_\exists =\left\langle \A_1^{\prime\prime} , \ldots ,\A_{\kappa}^{\prime\prime} \right\rangle$,
\linebreak
$\bb_\forall =\left\langle \bb_1^\prime , \ldots ,\bb_{\kappa}^\prime \right\rangle$, 
$\bb_\exists =\left\langle \bb_1^{\prime\prime} , \ldots ,\bb_{\kappa}^{\prime\prime} \right\rangle$, 
где $\kappa$ -- число AE-блоков в $Q(\A,\bb, \cal A, \beta )$ и для каждого AE-блока $B_s$ 
две интервальные матрицы 
$\A_s^\prime =(\ab_{sij}^\prime)\; , \A_s^{\prime\prime}=(\ab_{sij}^{\prime\prime})$ и два 
интервальных вектора  $\bb_s^\prime =(\bb_{si}^\prime )\; , \bb_s^{\prime\prime} =(\bb_{si}^{\prime\prime} )$ 
определяются следующим образом:

1) если либо $s=\overline{2,\kappa -1}$, либо $s=1$ и $B_1$ удовлетворяет условию 1),  либо $s=\kappa$ и $B_{\kappa}$ удовлетворяет условию 1), то полагаем 
$$
\ab_{sij}^{\prime}=
\left\{
\begin{array}{cl}
\ab_{ij}, & \mbox{если}\; a_{ij} \in \left\{ u_{i_{s}}, \ldots ,  u_{i_{s-1}+1} \right\} \& {\cal A}_{ij}=\forall, \\
{\bf 0}, & \mbox{в противном случае},
\end{array}
\right. 
$$
$$
\ab_{sij}^{\prime\prime}=
\left\{
\begin{array}{cl}
\ab_{ij}, & \mbox{если}\; a_{ij} \in \left\{ u_{i_{s}}, \ldots ,  u_{i_{s-1}+1} \right\} \& {\cal A}_{ij}=\exists, \\
{\bf 0}, & \mbox{в противном случае},
\end{array}
\right. 
$$
$$
\bb_{si}^{\prime}=
\left\{
\begin{array}{cl}
\bb_{i}, & \mbox{если}\; b_{i} \in \left\{ u_{i_{s}}, \ldots ,  u_{i_{s-1}+1} \right\} \& \beta_{i}=\forall, \\
{\bf 0}, & \mbox{в противном случае},
\end{array}
\right. 
$$
$$
\bb_{si}^{\prime\prime}=
\left\{
\begin{array}{cl}
\bb_{i}, & \mbox{если}\; b_{i} \in \left\{ u_{i_{s}}, \ldots ,  u_{i_{s-1}+1} \right\} \& \beta_{i}=\exists, \\
{\bf 0}, & \mbox{в противном случае},
\end{array}
\right. 
$$

2) если $B_1$ состоит только из кванторов всеобщности, то полагаем
$$
\ab_{1ij}^{\prime}=
\left\{
\begin{array}{cl}
\ab_{ij}, & \mbox{если}\; a_{ij} \in \left\{ u_{i_{1}}, \ldots ,  u_{1} \right\} \\
{\bf 0}, & \mbox{в противном случае}
\end{array}
\right. ,\,\;\;\; \A_1^{\prime\prime}={\bf 0},
$$
$$
\bb_{1i}^{\prime}=
\left\{
\begin{array}{cl}
\bb_{i}, & \mbox{если}\; b_{i} \in \left\{ u_{i_{1}}, \ldots ,  u_{1} \right\} \\
{\bf 0}, & \mbox{в противном случае}
\end{array}
\right. ,\,\;\; \; \bb_{1}^{\prime\prime}={\bf 0},
$$

3) если $B_{\kappa}$ состоит только из кванторов существования, то полагаем 
$$
\A_{\kappa}^{\prime}={\bf 0}, \;\;\;\;
\ab_{\kappa ij}^{\prime\prime}=
\left\{
\begin{array}{cl}
\ab_{ij}, & \mbox{если}\; a_{ij} \in \left\{ u_{\mu}, \ldots ,  u_{i_{\kappa-1}+1} \right\}, \\
{\bf 0}, & \mbox{в противном случае},
\end{array}
\right.  
$$
$$
\bb_{\kappa}^{\prime}={\bf 0}, \;\;\;\;
\bb_{\kappa i}^{\prime\prime}=
\left\{
\begin{array}{cl}
\bb_{i}, & \mbox{если}\; b_{i} \in \left\{ u_{\mu}, \ldots ,  u_{i_{\kappa-1}+1} \right\}, \\
{\bf 0}, & \mbox{в противном случае}.
\end{array}\right.   
$$

Нетрудно показать, что для так определенных матриц и векторов выполняются равенства 
$$
\A = 
\sum\limits_{s=1}^{\kappa}\left( \A_s^{\prime} + \A_s^{\prime\prime}\right), \;\;\;\; 
\bb =\sum\limits_{s=1}^{\kappa}\left( \bb_s^{\prime} + \bb_s^{\prime\prime} \right),
\eqno (6)
$$
матрицы $\A_1^\prime , \ldots ,\A_{\kappa}^\prime,\A_1^{\prime\prime} , \ldots ,\A_{\kappa}^{\prime\prime}$ образуют 
дизъюнктное разбиение матрицы $\A$, то есть при фиксированных $i, \, j$ среди интервалов 
$\ab_{1ij}^\prime , \ldots ,\ab_{\kappa ij}^\prime,\ab_{1ij}^{\prime\prime} , \ldots ,\ab_{\kappa ij}^{\prime\prime}$ 
не более одного ненулевого, а векторы 
$\bb_1^\prime , \ldots ,\bb_{\kappa}^\prime,\bb_1^{\prime\prime} , \ldots ,\bb_{\kappa}^{\prime\prime}$ образуют 
дизъюнктное разбиение вектора $\bb$, 
то есть при фиксированном $i$ среди интервалов 
$\bb_{1i}^\prime , \ldots ,\bb_{\kappa i}^\prime,\bb_{1i}^{\prime\prime} , \ldots ,\bb_{\kappa i}^{\prime\prime}$ 
не более одного ненулевого.

Заметим также, что кванторная приставка $Q(\A,\bb, \cal A, \beta )$ практически однозначно (с точностью до нулевых 
элементов $\A$ и $\bb$, которые не влияют на множество решений) определяется по кортежам 
$\A_\forall,\, \A_\exists,\, \bb_\forall$ и $ \bb_\exists$. Поэтому мы можем переформулировать определение 3 следующим, эквивалетным, образом.

{\bf Определение $3^\prime$}. Пусть для уравнения (1) заданы кортежи 
$\A_\forall,\, \A_\exists,\, \bb_\forall,\, \bb_\exists$, образующие дизъюнктные разбиения матрицы $\A$ и вектора $\bb$ 
соответственно. Ин\-тер\-валь\-но-кванторной системой линейных уравнений будем называть формулу
$$
\begin{array}{c}
\Phi_{\kappa} (\A_\forall , \A_\exists , \bb_\forall , \bb_\exists ,   x ) \equiv
(\forall A_{\kappa}^\prime \in \A_{\kappa}^\prime)(\forall b_{\kappa}^\prime \in \bb_{\kappa}^\prime) 
(\exists A_{\kappa}^{\prime\prime} \in \A_{\kappa}^{\prime\prime})(\exists b_{\kappa}^{\prime\prime} \in \bb_{\kappa}^{\prime\prime}) \ldots \\ \\
\ldots (\forall A_{1}^\prime \in \A_1^\prime) 
(\forall b_{1}^\prime \in \bb_{1}^\prime) 
(\exists A_{1}^{\prime\prime} \in \A_{1}^{\prime\prime})(\exists b_{1}^{\prime\prime} \in \bb_{1}^{\prime\prime}) 
\left(\sum\limits_{s=1}^{\kappa}\left( A_s^\prime + A_s^{\prime\prime} \right) x = 
\sum\limits_{s=1}^{\kappa}\left( b_s^\prime + b_s^{\prime\prime} \right)\right)
\end{array}
\eqno (7)
$$
и, соответственно, решением будет вектор $x \in \R^n$, для которого формула принимает
зна\-че\-ние “истина”, а множество решений (по определению)  (interval-quantifier solution set) 
$$
\varXi_{iq}^{\kappa}(\A, \bb, \A_\forall , \A_\exists , \bb_\forall , \bb_\exists )=\{x\in \R^n \; | \; 
\Phi_{\kappa} (\A_\forall , \A_\exists , \bb_\forall , \bb_\exists , x )=\mbox{"'истина''} \} .
\eqno (8)
$$

Представление интервально-кванторной системы линейных уравнений в виде оп\-ре\-де\-ле\-ния $3^\prime$, так же как и в случае 
AE-решений, наводит на мысль перейти к более общим уравнениям, чем (1) и (4). 

Будем считать, что заданы два кортежа интервальных матриц   
$\A_\forall =\left\langle \A_1^\prime , \ldots ,\A_{\kappa}^\prime \right\rangle$, 
$\A_\exists =\left\langle \A_1^{\prime\prime} , \ldots ,\A_{\kappa}^{\prime\prime} \right\rangle$
и два кортежа интервальных векторов
$\bb_\forall =\left\langle \bb_1^\prime , \ldots ,\bb_{\kappa}^\prime \right\rangle$, 
\linebreak
$\bb_\exists =\left\langle \bb_1^{\prime\prime} , \ldots ,\bb_{\kappa}^{\prime\prime} \right\rangle$. 
Рассмотрим следующее уравнение: 
$$
\left(\sum\limits_{i=1}^{\kappa}\A_i^\prime + \sum\limits_{i=1}^{\kappa}\A_i^{\prime\prime}\right)  x =
\sum\limits_{i=1}^{\kappa}\bb_i^\prime + \sum\limits_{i=1}^{\kappa}\bb_i^{\prime\prime}.
\eqno (9)
$$
При этом матрицы и векторы -- произвольные (не связанные соотношениями дизъюнктности).  

{\bf Определение 5}. Множеством решений уравнения (9) называется множество
$$
\begin{array}{l}
\varXi_{iq}^{\kappa}(\A_\forall , \A_\exists , \bb_\forall , \bb_\exists )=\{ x\in \R^n \; | \;
(\forall A_{\kappa}^\prime \in \A_{\kappa}^\prime)(\forall b_{\kappa}^\prime \in \bb_{\kappa}^\prime) 
(\exists A_{\kappa}^{\prime\prime} \in \A_{\kappa}^{\prime\prime})(\exists b_{\kappa}^{\prime\prime} 
\in \bb_{\kappa}^{\prime\prime}) \ldots \\ \\
\ldots (\forall A_{1}^\prime \in \A_1^\prime ) 
(\forall b_{1}^\prime \in \bb_{1}^\prime) 
(\exists A_{1}^{\prime\prime} \in \A_{1}^{\prime\prime})(\exists b_{1}^{\prime\prime} \in \bb_{1}^{\prime\prime}) 
(\sum\limits_{s=1}^{\kappa}( A_s^\prime + A_s^{\prime\prime} ) x = 
\sum\limits_{s=1}^{\kappa}( b_s^\prime + b_s^{\prime\prime} ) ) \} ,
\end{array}
\eqno (10)
$$
или, используя обозначение (7),
$$
\varXi_{iq}^{\kappa}(\A_\forall , \A_\exists , \bb_\forall , \bb_\exists )=\{x\in \R^n \; | \; 
\Phi_{\kappa} (\A_\forall , \A_\exists , \bb_\forall , \bb_\exists , x )=\mbox{"'истина''} \} .
$$

Таким образом уравнение (9) является некоторым обобщением понятия ин\-тер\-валь\-но-кван\-тор\-ной системы линейных уравнений,  
вводимых определением 3 (или $3^\prime$) (отсутствуют условия дизъюнктности). Поэтому уравнение (9) будем называть 
обобщенным интервально-кванторным линейным уравнением. 

Основная задача данной работы -- получить  бескванторное описание именно множества 
$\varXi_{iq}^{\kappa}(\A_\forall , \A_\exists , \bb_\forall , \bb_\exists )$ и, как следствие,  
бескванторное описание множества 
\linebreak
$\varXi_{iq}^{\kappa}(\A, \bb, \A_\forall , \A_\exists , \bb_\forall , \bb_\exists )$.

Отметим также, что при $\kappa = 1$ уравнение (9) совпадает с уравнением (4).

\begin{center}
{\bf 2. Обобщенные интервально-кванторные линейные  уравнения \\
в упорядоченных векторных пространствах}
\end{center}

Заметим, что понятие обобщенного интервально-кванторного линейного  уравнения в виде уравнения (9), 
в отличие от понятия ин\-тер\-валь\-но-кван\-тор\-ной системы линейных уравнений,  
вводимых определением 3 (или $3^\prime$),
легко переносится на произвольные упорядоченные векторные пространства, при этом будем придерживаться 
терминологии из \cite{AK,V}. 

Пусть $\kb{E}$ -- упорядоченное векторное  пространство над $\kb{R}^1$ и $\kb{E}_+$ -- конус положительных элементов 
в $\kb{E}$. По аналогии с 
$\kb{IR}$ обозначим $\kb{IE}$  множество интервалов в $\kb{E}$, т.е. если $\un{a},\ov{a}\in \kb{E}, \;\; \un{a}\leq \ov{a}$, то 
$\ab = [\un{a},\ov{a}]=\{ x\in \kb{E} \;|\; \un{a}\leq x \leq\ov{a} \}\in \kb{IE} $. 
При этом в дальнейшем, как обычно, будем отождествлять одноэлементные интервалы вида $[a,a]$ с элементом  $a\in \kb{E}$ и считать $\kb{E}\subseteq \kb{IE}$.

Если $\kb{F}$ -- также упорядоченное векторное  пространство и 
$\kb{L}(\kb{E},\kb{F})$ -- пространство линейных операторов из $\kb{E}$ в $\kb{F}$, то будем считать его предупорядоченным с помощью конуса положительных операторов и обозначать $\kb{IL(E,F)}$ -- множество интервалов в нем. Известно, что это предупорядочение 
будет упорядочением,  если конус $\kb{E}_+$ -- воспроизводящий \cite[стр. 108]{AK}.

Формально заменяя в формуле (9) интервальные матрицы  на интервальные операторы из $\kb{IL(E,F)}$, 
а интервальные векторы -- на интервалы из 
$\kb{IF}$, получаем следующее определение обобщеного интервально-кванторного линейного  уравнения в упорядоченных векторных пространствах.

Пусть заданы два кортежа интервальных операторов   
$\A_\forall =\left\langle \A_1^\prime , \ldots ,\A_{\kappa}^\prime \right\rangle$,
\linebreak 
$\A_\exists =\left\langle \A_1^{\prime\prime} , \ldots ,\A_{\kappa}^{\prime\prime} \right\rangle$
и два кортежа интервалов
$\bb_\forall =\left\langle \bb_1^\prime , \ldots ,\bb_{\kappa}^\prime \right\rangle$, 
\linebreak
$\bb_\exists =\left\langle \bb_1^{\prime\prime} , \ldots ,\bb_{\kappa}^{\prime\prime} \right\rangle$. 
Рассмотрим следующее уравнение: 
$$
\left(\sum\limits_{i=1}^{\kappa}\A_i^\prime + \sum\limits_{i=1}^{\kappa}\A_i^{\prime\prime}\right)  x =
\sum\limits_{i=1}^{\kappa}\bb_i^\prime + \sum\limits_{i=1}^{\kappa}\bb_i^{\prime\prime}.
\eqno (9a)
$$

{\bf Определение 5a}. Множеством решений уравнения (9a) называется множество
$$
\begin{array}{l}
\varXi_{iq}^{\kappa}(\A_\forall , \A_\exists , \bb_\forall , \bb_\exists )=\{ x\in \R^n \; | \;
(\forall A_{\kappa}^\prime \in \A_{\kappa}^\prime)(\forall b_{\kappa}^\prime \in \bb_{\kappa}^\prime) 
(\exists A_{\kappa}^{\prime\prime} \in \A_{\kappa}^{\prime\prime})(\exists b_{\kappa}^{\prime\prime} 
\in \bb_{\kappa}^{\prime\prime}) \ldots \\ \\
\ldots (\forall A_{1}^\prime \in \A_1^\prime ) 
(\forall b_{1}^\prime \in \bb_{1}^\prime) 
(\exists A_{1}^{\prime\prime} \in \A_{1}^{\prime\prime})(\exists b_{1}^{\prime\prime} \in \bb_{1}^{\prime\prime}) 
(\sum\limits_{s=1}^{\kappa}( A_s^\prime + A_s^{\prime\prime} ) x = 
\sum\limits_{s=1}^{\kappa}( b_s^\prime + b_s^{\prime\prime} ) ) \} ,
\end{array}
\eqno (10a)
$$
или, используя обозначение (7), 
$$
\varXi_{iq}^{\kappa}(\A_\forall , \A_\exists , \bb_\forall , \bb_\exists )=\{x\in \R^n \; | \; 
\Phi_{\kappa} (\A_\forall , \A_\exists , \bb_\forall , \bb_\exists , x )=\mbox{"'истина''} \} .
$$

\begin{center}
{\bf 3. Бескванторное описание множества решений обобщенного \\ интервально-кванторного  
линейного  уравнения в интервальной арифметике}
\end{center}

Для исследования интервальных уравнений в упорядоченных векторных пространствах нам понадобятся два 
утверждения: 1) что сумма двух интервалов -- интервал;
2) что если $\A \in \kb{IL(E,F)}$ -- интервальный 
оператор и $x \in \kb{E}$, то $\A x$ -- интервал в $\kb{F}$, которые, вообще говоря, не верны для произвольных упорядоченных векторных пространств (соответствующие примеры имеются в \cite{lak1, lak2}). Поэтому в данном разделе мы будем рассматривать их как дополнительные предположения, а в следующем приведем некоторые достаточные для них условия. 
Отметим, что хотя в общем случае $\ab + \bb$ может не быть интервалом в $\kb{F}$, но если хотя бы один из интервалов 
$\ab$ или  $\bb$ одноэлементный, то $\ab + \bb$, очевидно, интервал.  
При этом, как обычно, сумма интервалов $\ab \in\kb{IF}$ и  $\bb\in\kb{IF}$  определяется как суммы всех элементов из $\ab$ и $\bb$, т.е. 
$\ab + \bb = \{ x+y \;|\; x\in \ab, y\in\bb \}$, и $\A x=\{ Ax \;|\; A\in \A \}$. В случае, когда $\kb{E}$ и $\kb{F}$ -- конечномерные пространства с обычным (покоординатным) упорядочением, множество $\A x$ легко вычисляется по формулам умножения интервальной матрицы на вектор \cite{sarb}.

{\bf Определение 3.} Будем говорить, что в упорядоченном векторном пространстве $\kb{F}$ выполняется теорема о сумме интервалов, если для любых $\ab =[\un{a} ,\ov{a}  ] \in\kb{IF}$ и  
\linebreak
$\bb=[\un{b} ,\ov{b}] \in\kb{IF}$ их сумма $\ab + \bb $ тоже интервал, при этом очевидно, что должно выполнятся равенство  $\ab + \bb =[\un{a}+\un{b} ,\ov{a}+\ov{b}]$.

Известно, что это так для векторных решеток  \cite[теорема 3(2.1), стр. 114]{AK} и для упорядоченных векторных
пространств,  обладающих свойством $C$ \cite{lak1, lak2}, определяемых в следующем разделе. 

Далее нам также понадобятся следующие утверждения. 

{\bf Лемма 1.} Пусть  $\kb{F}$ -- упорядоченное векторное пространство. Тогда для него следующие утверждения зквивилентны:
\newline
\indent
1) в $\kb{F}$ выполняется теорема о сумме интервалов; 
\newline
\indent
2) для любых $\ab =[\un{a} ,\ov{a}  ] \in\kb{IF}$ и  $\bb=[\un{b} ,\ov{b}] \in\kb{IF}$ выполняется 
$\ab \cap \bb \neq \emptyset$ тогда и только тогда, когда $\un{a}\leq \ov{b}$ и $\un{b}\leq \ov{a} $;
\newline
\indent
3) если $a_1,a_2,b_1,b_2\in \kb{F}$ такие, что $a_i\leq b_j$, для всех $i,j=1,2$, то существует $c\in \kb{F}$  такое, что 
$a_i\leq c \leq b_j$ для всех $i,j=1,2$.

{\bf Доказательство.} $1)\Rightarrow 2).$ Берем $\ab =[\un{a} ,\ov{a}  ] \in\kb{IF}$ и  $\bb=[\un{b} ,\ov{b}] \in\kb{IF}$.
Если $\ab \cap \bb \neq \emptyset$ то найдется $x \in\ab \cap \bb $. Поэтому $\un{a}\leq x\leq \ov{b}$ и 
$\un{b}\leq x\leq \ov{a} $.
\newline
\indent
Обратно, пусть $\un{a}\leq \ov{b}$ и $\un{b}\leq \ov{a} $. По условию 1) $\ab - \bb $ -- интервал,  
$\ab - \bb = [\un{a}-\ov{b} ,\ov{a}-\un{b}]$ и так как $\un{a}- \ov{b}\leq 0 \leq \ov{a}-\un{b}$, то $0\in\ab - \bb $ . Но тогда найдутся $x_1\in \ab $ и $x_2\in \bb$ такие, что $0=x_1-x_2$. Следовательно $x_1=x_2\in \ab \cap \bb \neq \emptyset$.
\newline
\indent
$2)\Rightarrow 1).$ Пусть выполняется условие 2). Берем $\ab =[\un{a} ,\ov{a}  ] \in\kb{IF}$,   
$\bb=[\un{b} ,\ov{b}] \in\kb{IF}$ и покажем, что $\ab + \bb $ -- интервал. Так как всегда 
$\ab + \bb \subseteq [\un{a}+\un{b} ,\ov{a}+\ov{b}]$, то достаточно показать обратное включение, то есть
$[\un{a}+\un{b} ,\ov{a}+\ov{b}]\subseteq \ab + \bb $. 
\newline
\indent
Пусть $x\in [\un{a}+\un{b} ,\ov{a}+\ov{b}]$. Рассмотрим интервалы $[x-\ov{a}, x-\un{a}]$ и $\bb=[\un{b} ,\ov{b}]$. 
Для них выполнены условия $x-\ov{a}\leq\ov{b}$ и $\un{b}\leq x-\un{a}$. Поэтому из условия 2) получаем, что 
$[x-\ov{a}, x-\un{a}]\cap [\un{b} ,\ov{b}] \neq \emptyset$. Пусть $x_1\in [x-\ov{a}, x-\un{a}]\cap [\un{b} ,\ov{b}]$. 
Обозначая $x_2 =x -x_1$, получаем, что $\un{a}\leq x_2 =x -x_1\leq \ov{a}$, то есть $x_2 \in \ab =[\un{a} ,\ov{a} ]$. 
Следовательно, $x=x_2 +x_1 \in \ab + \bb$.
\newline
\indent
Эквивалентность условий 2) и 3) легко получается, если положить $a_1 =\un{a},\;a_2 =\un{b},$ 
\linebreak 
$b_1 =\ov{a}, \;b_2 =\ov{b}$. Лемма доказана.

Будем использовать в дальнейшем следующее понятие ширины интервала
\linebreak 
 $\ab =[\un{a} ,\ov{a}  ] \in\kb{IF}$, обозначаемое 
${\rm wid}(\ab) = \ov{a} - \un{a}$. Если в $\kb{F}$ выполняется теорема о сумме интервалов, то для функции 
${\rm wid}$, очевидно, выполняются те же свойства, что и в $\kb{IR}^n$ \cite{sarb, Neum}, в частности, 
${\rm wid}(\ab + \bb) = {\rm wid}(\ab)+{\rm wid}(\bb)$ и  ${\rm wid}(\ab + b) = {\rm wid}(\ab)$, если $b\in\kb{F}$.

{\bf Лемма 2.} Пусть для $\kb{F}$ выполняется теорема о сумме интервалов и $\ab ,\bb ,\cb \in\kb{IF}$. 
\linebreak 
Тогда в $\kb{IF}$ следующие условия эквивалентны: 
\newline
\indent
1) $(\exists c\in\cb) \;(\ab \subseteq \bb +c) $; 
\newline
\indent
2) $(\ab \subseteq \bb +\cb )\; \& \;{\rm wid}(\ab)\leq {\rm wid}(\bb)$.

{\bf Доказательство.} Пусть выполняется условие $(\exists c\in\cb) \;(\ab \subseteq \bb +c) $. Берем $c\in\cb$ такое, что 
 $\ab \subseteq \bb +c$. Тогда очевидно, что ${\rm wid}(\ab)\leq {\rm wid}(\bb+c) = {\rm wid}(\bb)$ и 
\linebreak 
$\ab \subseteq \bb +c\subseteq \bb +\cb$. 
\newline
\indent
Покажем обратное. Пусть $(\ab \subseteq \bb +\cb )\; \& \;{\rm wid}(\ab)\leq {\rm wid}(\bb)$.
\newline
\indent
Из неравенства ${\rm wid}(\ab)\leq {\rm wid}(\bb)$ получаем, что $\ov a -\un a \leq \ov b -\un b$, или 
$\ov a -\ov b \leq \un a -\un b$. Поэтому можно определить интервал $[\ov a -\ov b , \un a -\un b]$ . 
Из включения 
$\ab\subseteq \bb + {\bf c}$ и теоремы о сумме интервалов получаем, что выполняются неравенства 
$\un b + \un c\leq \un a $ и $\ov a \leq \ov b +\ov c$, или $ \un c\leq \un a -\un b $ и $\ov a -\ov b \leq  \ov c$. 
Но тогда для  интервалов $[\ov a -\ov b , \un a -\un b]$ и ${\bf c}=[\un c ,\ov c]$ выполняется условие 2) леммы 1 и, следовательно, они имеют не пустое пересечение. 
Поэтому найдется $c_0 \in {\bf c}$ такое, что $\ov a -\ov b \leq c_0 \leq \un a -\un b$, или $\un b + c_0 \leq \un a$,  
$\ov a  \leq \ov b +c_0$. Отсюда  получаем, что $\ab\subseteq \bb + c_0$. Лемма доказана.

Бескванторное описание множества $\varXi_{iq}^{\kappa}(\A_\forall , \A_\exists , \bb_\forall , \bb_\exists )$ содержится в следующей теореме.

{\bf Теорема 3.} Пусть для упорядоченных векторных пространств $\kb{E}$ и $\kb{F}$  
заданы два кортежа интервальных операторов 
$\A_\forall =\left\langle \A_1^\prime , \ldots ,\A_{\kappa}^\prime \right\rangle$, 
$\A_\exists =\left\langle \A_1^{\prime\prime} , \ldots ,\A_{\kappa}^{\prime\prime} \right\rangle$, где
\linebreak
$\A_i^\prime ,\A_i^{\prime\prime} \in \kb{IL(E,F)}$, 
два кортежа интервалов 
$\bb_\forall =\left\langle \bb_1^\prime , \ldots ,\bb_{\kappa}^\prime \right\rangle$, 
$\bb_\exists =\left\langle \bb_1^{\prime\prime} , \ldots ,\bb_{\kappa}^{\prime\prime} \right\rangle$, где 
$\bb_i^\prime ,\bb_i^{\prime\prime} \in \kb{IF}$, ($i=\overline{1,\kappa}$) 
и выполняются следующие условия:
\newline
\indent
1) для любого $x \in\kb{E}$ и всех $i=\ov{1,\kappa}\;\;$ $\A_i^\prime x, \; \A_i^{\prime\prime} x$ -- интервалы в $\kb{F}$;
\newline
\indent
2) для $\kb{F}$ выполняется теорема о сумме интервалов.
\newline
\indent 
Тогда $x\in \varXi_{iq}^{\kappa}(\A_\forall , \A_\exists , \bb_\forall , \bb_\exists )$ если и только если выполняются 
следующие условия:
$$
\sum\limits_{i=1}^{\kappa}\left(\A_i^\prime x - \bb_i^\prime\right) 
\subseteq 
\sum\limits_{i=1}^{\kappa}\left(\bb_i^{\prime\prime} - \A_i^{\prime\prime}  x \right)  ,
\eqno (11)
$$
$$
\sum\limits_{i=1}^{l}{\rm wid}(\A_i^\prime x - \bb_i^\prime) 
\leq 
\sum\limits_{i=1}^{l}{\rm wid}(\bb_i^{\prime\prime} - \A_i^{\prime\prime}  x) , 
\;\;\;\; l= \ov {1, \kappa -1}.
\eqno (12)
$$

{\bf Доказательство.} Доказательство проведем индукцией по $\kappa$. 
\newline
\indent
База индукции. При $\kappa =1$  множество $\varXi_{iq}^{\kappa}(\A_\forall , \A_\exists , \bb_\forall , \bb_\exists )$ имеет вид
$$
\begin{array}{l}
\varXi_{iq}^{1}(\A_\forall , \A_\exists , \bb_\forall , \bb_\exists )= 
\\ \\
=\{x\in \kb{E} \; | \; 
(\forall A_{1}^\prime \in \A_1^\prime) 
(\forall b_{1}^\prime \in \bb_{1}^\prime) 
(\exists A_{1}^{\prime\prime} \in \A_{1}^{\prime\prime})(\exists b_{1}^{\prime\prime} \in \bb_{1}^{\prime\prime}) 
\left( A_1^\prime + A_1^{\prime\prime} \right) x = 
b_1^\prime + b_1^{\prime\prime}\}.
\end{array} 
$$

Поэтому $x\in \varXi_{iq}^{1}(\A_\forall , \A_\exists , \bb_\forall , \bb_\exists )$ означает, что при любых 
$A_{1}^\prime \in \A_1^\prime,\; b_{1}^\prime \in \bb_{1}^\prime$ уравнение 
$\left( A_1^\prime + A_1^{\prime\prime} \right) x = b_1^\prime + b_1^{\prime\prime}$ или, что то же самое, 
$ A_1^\prime x - b_1^\prime  =   b_1^{\prime\prime} -A_1^{\prime\prime} x$ выполняется при некоторых 
$A_{1}^{\prime\prime} \in \A_{1}^{\prime\prime}$ и $b_{1}^{\prime\prime} \in \bb_{1}^{\prime\prime}$.
Следовательно, при любых 
$A_{1}^\prime \in \A_1^\prime,\; b_{1}^\prime \in \bb_{1}^\prime$
$$
 A_1^\prime x - b_1^\prime  \in \{  b_1^{\prime\prime} -A_1^{\prime\prime} x\;|\; 
A_{1}^{\prime\prime} \in \A_{1}^{\prime\prime}, \;b_{1}^{\prime\prime} \in \bb_{1}^{\prime\prime}\},
$$
а значит, и 
$$
\{ A_1^\prime x - b_1^\prime \;|\;A_{1}^\prime \in \A_1^\prime,\; b_{1}^\prime \in \bb_{1}^\prime\}  
\subseteq \{  b_1^{\prime\prime} -A_1^{\prime\prime} x\;|\; 
A_{1}^{\prime\prime} \in \A_{1}^{\prime\prime}, \;b_{1}^{\prime\prime} \in \bb_{1}^{\prime\prime}\},
$$ 
то есть 
$\left(\A_1^\prime x - \bb_1^\prime\right) 
\subseteq 
\left(\bb_1^{\prime\prime} - \A_1^{\prime\prime}  x \right) $ (при $\kappa =1$ условие (12) отсутствует).

Пусть теорема верна для $(\kappa -1)$ ($\kappa\geq 2$) и покажем, что она верна и для $\kappa$.

Итак, пусть заданы кортежи 
$\A_\forall =\left\langle \A_1^\prime , \ldots ,\A_{\kappa}^\prime \right\rangle$, 
$\A_\exists =\left\langle \A_1^{\prime\prime} , \ldots ,\A_{\kappa}^{\prime\prime} \right\rangle$, 
\linebreak 
$\bb_\forall =\left\langle \bb_1^\prime , \ldots ,\bb_{\kappa}^\prime \right\rangle$, 
$\bb_\exists =\left\langle \bb_1^{\prime\prime} , \ldots ,\bb_{\kappa}^{\prime\prime} \right\rangle$.
Возьмем любые $A_{\kappa}^\prime, A_{\kappa}^{\prime\prime}\in \kb{L(E,F)}$ и $b_{\kappa}^\prime, b_{\kappa}^{\prime\prime}\in \kb{F}$ и образуем кортежи
$\tilde{\A}_\forall =\left\langle \tilde{\A}_1^\prime , \ldots ,\tilde{\A}_{\kappa -1}^\prime \right\rangle$, 
$\tilde{\A}_\exists =\left\langle \tilde{\A}_1^{\prime\prime} , \ldots ,\tilde{\A}_{\kappa -1}^{\prime\prime} \right\rangle$, 
$\tilde{\bb}_\forall =\left\langle \tilde{\bb}_1^\prime , \ldots ,\tilde{\bb}_{\kappa -1}^\prime \right\rangle$, 
$\tilde{\bb}_\exists =\left\langle \tilde{\bb}_1^{\prime\prime} , \ldots ,\tilde{\bb}_{\kappa -1}^{\prime\prime} \right\rangle$
длины $(\kappa -1)$, полагая
$$
\tilde{\A}_{\kappa -1}^\prime =\A_{\kappa -1}^\prime +A_{\kappa}^\prime , \;
\tilde{\A}_{\kappa -1}^{\prime\prime} = \A_{\kappa -1}^{\prime\prime} +A_{\kappa}^{\prime\prime} , \; 
\tilde{\bb}_{\kappa -1}^\prime =\bb_{\kappa -1}^\prime +b_{\kappa}^\prime , \;
\tilde{\bb}_{\kappa -1}^{\prime\prime} = \bb_{\kappa -1}^{\prime\prime} +b_{\kappa}^{\prime\prime} , 
$$
и 
$
\tilde{\A}_{i}^\prime =\A_{i}^\prime , \;
\tilde{\A}_{i}^{\prime\prime} = \A_{i}^{\prime\prime} , \; 
\tilde{\bb}_{i}^\prime =\bb_{i}^\prime , \;
\tilde{\bb}_{i}^{\prime\prime} = \bb_{i}^{\prime\prime}  
$
для $i=\ov{1,\kappa -2}$. Очевидно, что при этом 
условие 2) теоремы 3 будет выполняться и для новых кортежей при любых
\linebreak  
$A_{\kappa}^\prime, A_{\kappa}^{\prime\prime}\in \kb{L(E,F)}$ и $b_{\kappa}^\prime, b_{\kappa}^{\prime\prime}\in \kb{F}$.

Заметим также, что для формул, определяющих множества  $\varXi_{iq}^{\kappa}(\A_\forall , \A_\exists , \bb_\forall , \bb_\exists )$ и
\linebreak 
$\varXi_{iq}^{\kappa -1}(\tilde{\A}_\forall , \tilde{\A}_\exists , \tilde{\bb}_\forall , \tilde{\bb}_\exists )$, 
верно следующее утверждение:
$$
\begin{array}{l}
\Phi_{\kappa} (\A_\forall , \A_\exists , \bb_\forall , \bb_\exists ,   x ) \Leftrightarrow \\ \\
\Leftrightarrow (\forall A_{\kappa}^\prime \in \A_{\kappa}^\prime)(\forall b_{\kappa}^\prime \in \bb_{\kappa}^\prime) 
(\exists A_{\kappa}^{\prime\prime} \in \A_{\kappa}^{\prime\prime})(\exists b_{\kappa}^{\prime\prime} \in \bb_{\kappa}^{\prime\prime}) \,
\Phi_{\kappa -1} (\tilde{\A}_\forall , \tilde{\A}_\exists , \tilde{\bb}_\forall , \tilde{\bb}_\exists ,   x ), 
\end{array}
\eqno (13)
$$ 
где $\Leftrightarrow$ -- логическая эквивалентность формул.

По индуктивному предположению получаем, что 
$x\in\varXi_{iq}^{\kappa -1}(\tilde{\A}_\forall , \tilde{\A}_\exists , \tilde{\bb}_\forall , \tilde{\bb}_\exists )$ 
тогда и только тогда, когда выполняются условия 
$$
\sum\limits_{i=1}^{\kappa -1}\left(\tilde{\A}_i^\prime x - \tilde{\bb}_i^\prime\right) 
\subseteq 
\sum\limits_{i=1}^{\kappa -1}\left(\tilde{\bb}_i^{\prime\prime} - \tilde{\A}_i^{\prime\prime}  x \right)  ,
\eqno (11^\prime)
$$
$$
\sum\limits_{i=1}^{l}{\rm wid}(\tilde{\A}_i^\prime x - \tilde{\bb}_i^\prime) 
\leq 
\sum\limits_{i=1}^{l}{\rm wid}(\tilde{\bb}_i^{\prime\prime} - \tilde{\A}_i^{\prime\prime}  x) , 
\;\;\;\; l= \ov {1, \kappa -2},
\eqno (12^\prime)
$$
которые после подстановки в них выражений для 
$
\tilde{\A}_{i}^\prime  , \;
\tilde{\A}_{i}^{\prime\prime}  , \; 
\tilde{\bb}_{i}^\prime  , \;
\tilde{\bb}_{i}^{\prime\prime} , 
$
 $i=\ov{1,\kappa -1}$ и простых преобразований принимают вид
$$
\sum\limits_{i=1}^{\kappa -1}\left(\A_i^\prime x - \bb_i^\prime\right) +A_{\kappa}^\prime x -b_{\kappa}^\prime
\subseteq 
\sum\limits_{i=1}^{\kappa -1}\left(\bb_i^{\prime\prime} - \A_i^{\prime\prime}  x \right) + b_{\kappa}^{\prime\prime} -A_{\kappa}^{\prime\prime}  x,
\eqno (11^{\prime\prime})
$$
$$
\sum\limits_{i=1}^{l}{\rm wid}(\A_i^\prime x - \bb_i^\prime) 
\leq 
\sum\limits_{i=1}^{l}{\rm wid}(\bb_i^{\prime\prime} - \A_i^{\prime\prime}  x) , 
\;\;\;\; l= \ov {1, \kappa -2}.
\eqno (12^{\prime\prime})
$$

Обозначим $\ab_0 = \sum\limits_{i=1}^{\kappa -1}\left(\A_i^\prime x - \bb_i^\prime\right)$, $\ab =\ab_0 +A_{\kappa}^\prime x -b_{\kappa}^\prime$, 
$\bb =\sum\limits_{i=1}^{\kappa -1}\left(\bb_i^{\prime\prime} - \A_i^{\prime\prime}  x \right)$,
\linebreak  
$\cb = \bb_{\kappa}^{\prime\prime} -\A_{\kappa}^{\prime\prime}  x$. Тогда из эквивалентности (13), учитывая, что условие 
$(12^{\prime\prime})$ не зависит от $A_{\kappa}^\prime,\; A_{\kappa}^{\prime\prime} ,\;b_{\kappa}^\prime ,\;b_{\kappa}^{\prime\prime}$, получаем

\noindent
$\Phi_{\kappa} (\A_\forall , \A_\exists , \bb_\forall , \bb_\exists ,  x ) \Leftrightarrow $

\noindent
$\Leftrightarrow ((\forall A_{\kappa}^\prime \in \A_{\kappa}^\prime)(\forall b_{\kappa}^\prime \in \bb_{\kappa}^\prime) 
(\exists A_{\kappa}^{\prime\prime} \in \A_{\kappa}^{\prime\prime})(\exists b_{\kappa}^{\prime\prime} \in \bb_{\kappa}^{\prime\prime}) 
(\ab_0 +A_{\kappa}^\prime x -b_{\kappa}^\prime \subseteq \bb + b_{\kappa}^{\prime\prime} -A_{\kappa}^{\prime\prime}  x)) \& (12^{\prime\prime}).$

Далее, так как очевидно, что 
$$
(\exists A_{\kappa}^{\prime\prime} \in \A_{\kappa}^{\prime\prime})(\exists b_{\kappa}^{\prime\prime} \in \bb_{\kappa}^{\prime\prime}) \,
(\ab_0 +A_{\kappa}^\prime x -b_{\kappa}^\prime \subseteq \bb + b_{\kappa}^{\prime\prime} -A_{\kappa}^{\prime\prime}  x) \Leftrightarrow 
(\exists c\in\cb) (\ab \subseteq \bb + c),
$$
то по лемме 2, учитывая, что ${\rm wid}(\ab)={\rm wid}(\ab_0)$, получаем 
$$
(\exists A_{\kappa}^{\prime\prime} \in \A_{\kappa}^{\prime\prime})(\exists b_{\kappa}^{\prime\prime} \in \bb_{\kappa}^{\prime\prime}) \,
(\ab_0 +A_{\kappa}^\prime x -b_{\kappa}^\prime \subseteq \bb + b_{\kappa}^{\prime\prime} -A_{\kappa}^{\prime\prime}  x) \Leftrightarrow 
 (\ab \subseteq \bb + \cb)\& {\rm wid}(\ab_0)\leq {\rm wid}(\bb).
$$

Поэтому
$$
\begin{array}{l}
\Phi_{\kappa} (\A_\forall , \A_\exists , \bb_\forall , \bb_\exists ,   x ) \Leftrightarrow 
\\ \\
\Leftrightarrow ((\forall A_{\kappa}^\prime \in \A_{\kappa}^\prime)(\forall b_{\kappa}^\prime \in \bb_{\kappa}^\prime) 
(\ab_0 +A_{\kappa}^\prime x -b_{\kappa}^\prime \subseteq \bb + \cb) \& ({\rm wid}(\ab_0)\leq {\rm wid}(\bb))
\& (12^{\prime\prime}).
\end{array} 
$$

Кроме того, понятно, что
$$
\begin{array}{l}
((\forall A_{\kappa}^\prime \in \A_{\kappa}^\prime)(\forall b_{\kappa}^\prime \in \bb_{\kappa}^\prime) 
(\ab_0 +A_{\kappa}^\prime x -b_{\kappa}^\prime \subseteq \bb + \cb) \Leftrightarrow 
(\ab_0 +\A_{\kappa}^\prime x -\bb_{\kappa}^\prime \subseteq \bb + \cb) \Leftrightarrow 
(11), \\ \\
({\rm wid}(\ab_0)\leq {\rm wid}(\bb))
\& (12^{\prime\prime}) \Leftrightarrow (12),
\end{array} 
$$
поэтому
$$
\Phi_{\kappa} (\A_\forall , \A_\exists , \bb_\forall , \bb_\exists ,   x ) \Leftrightarrow (11)\&(12).
$$
Теорема доказана.

{\bf Замечание 2.} Из доказательства видно, что при $\kappa =1$ условия 1) и 2) теоремы не нужны, а само  доказательство практически повторяет доказательство  теремы 1 из статьи С.П. Шарого \cite{sar4}.

{\bf Замечание 3.} Если $\kb{E}=\R^n$ и $\kb{F}=\R^m$ с покоординатными упорядочиваниями, то условия 1), 2) теоремы 1, очевидно, выполняются. При этом условие 2) тривиально, а условие 1) является следствием теоремы Оеттли-Прагера 
\cite{oetpr, Oet}.
Явную формулу для вычисления $\A x$ можно найти, например, в \cite[Proposition 2.27, стр. 62]{FNR} и в 
\cite[Предложение 2.2.4, стр. 94-95]{sarb}.
\newline
\indent
Кроме того, так как в этом случае упорядоченное постранство $\kb{IL(E,F)}$ изоморфно пространству $m\times n$-матриц с поэлементным упорядочением, то для него верна теорема о сумме интервалов и, следовательно, 
можно определить интервальные операторы
$\A^\prime = \sum\limits_{i=1}^{\kappa}\A_i^\prime $, $\A^{\prime\prime}= \sum\limits_{i=1}^{\kappa}\A_i^{\prime\prime}$ 
и интервальные векторы 
$\bb^\prime =\sum\limits_{i=1}^{\kappa} \bb_i^\prime$, $\bb^{\prime\prime}=\sum\limits_{i=1}^{\kappa}\bb_i^{\prime\prime}$, 
а условие (11) теоремы 3 записать в виде
$$
\A^\prime x - \bb^\prime \subseteq \bb^{\prime\prime} - \A^{\prime\prime} x .
\eqno (11a)
$$
Это условие в точности совпадает с необходимым и достаточным условием принадлежности вектора $x\in \R^n $ множеству 
$\varXi_{\forall\exists} (\A^{\prime}, \A^{\prime\prime}, \bb^{\prime}, \, \bb^{\prime\prime})$ из $2^\prime$.
Очевидно, что то же самое условие (11a) получится, если в формуле (10a) вынести вперед все кванторы всеобщности.

\begin{center}
{\bf 4. Теорема Оеттли-Прагера и  бескванторное описание множества \\
решений уравнения (9a) в виде неравенств с модулями }
\end{center}

Для того, чтобы получить описание множества решений через концы интервалов, необходимо иметь формулы, позволяющие найти интервал $\A x$ в $\kb{F}$ для $\A\in \kb{IL(E,F)}$. Поскольку для произвольных упорядоченных векторных пространств $\kb{E}$ и $\kb{F}$ множество $\A x$ может не быть интервалом, то нужны некоторые дополнительные условия на  $\kb{E}$ и $\kb{F}$. Оказалось, что для этого достаточно, чтобы в $\kb{E}$ и $\kb{F}$ существовало "`достаточно много"' мультипликаторов  
\cite[стр.219]{AK}. 

Оператор $\tau\in\kb{L(E,E)}$ такой, что $0_{\kb{E}}\leq\tau\leq I_{\kb{E}}$, где $0_{\kb{E}}$ и $I_{\kb{E}}$ -- нулевой  и тождественнный операторы в $\kb{E}$, соответственно, называется 
мультипликатором в $\kb{E}$. Множество всех мультипликаторов в $\kb{E}$ будем обозначать  ${\bf\Lambda} (\kb{E})$, то есть 
${\bf\Lambda} (\kb{E})=[0_{\kb{E}}, I_{\kb{E}}]\in\kb{IL(E,E)}$ \cite[стр.219]{AK}. Например, очевидно, что мультипликаторами в 
$\R^n$ будут линейные операторы, задаваемые диагональной матрицей, у которой на диагонали стоят числа из интервала $[0,1]$.

Всюду далее будем рассматривать  пространство линейных операторов из $\kb{E}$ в $\kb{F}$ только в случае, когда $\kb{E}$ является векторной решеткой \cite[стр.109]{AK}. Если $\kb{E}$ -- векторная решетка, то для любых $x, y\in \kb{E}$ через $x\vee y$ и $x\wedge y$ 
будем обозначать точные верхнюю и нижнюю границы элементов $x$ и $y$, $x^+ =x\vee 0$ -- положительная часть $x$, 
$x^- =-(x\wedge 0)$ -- отрицательная часть $x$ (при этом $x=x^+ - x^-$,  $x^+ \wedge x^- =0$) и $|x|=x^+ + x^- =x\vee (-x)$ -- модуль (или абсолютная величина)  элемента $x$. Так как в этом случае конус $\kb{E}_+$ воспроизводящий, то $\kb{L(E,F)}$ будет 
упорядочено с помощью конуса положительных операторов.

{\bf Определение 6.}
\newline
\indent 
1. Будем говорить, что упорядоченное векторное
пространство $\kb{F}$ обладает свойством $C$ (сжимаемости), если для любых $x, y\in \kb{F}_+$, $x\leq y$ существует мультипликатор $\tau\in {\bf\Lambda} (\kb{F})$, переводящий $y$ в $x$, т.е. такой, что $\tau y=x$. Класс упорядоченных векторных пространств, обладающих свойством $C$, обозначим через ${\cal K}(C)$. 
\newline
\indent
2. Будем говорить, что векторная решетка $\kb{E}$ обладает свойством $D$ (дизъюнктности), если для любых $x, y\in \kb{E}_+$, 
$x\wedge y=0$,
существует мультипликатор $\zeta\in {\bf\Lambda} (\kb{E})$ такой, что $\zeta x=x$ и $\zeta y=0$. Класс
векторных решеток, обладающих свойством $D$, обозначим через ${\cal K}(D)$.

Как уже упоминалось,  для векторных решеток выполняется теорема о сумме интервалов \cite[теорема 3(2.1), стр. 114]{AK}.
Покажем, что это так и для $\kb{F}\in {\cal K}(C)$. 

{\bf Лемма 3.} Если $\kb{F}\in {\cal K}(C)$, то для $\kb{F}$ выполняется теорема о сумме интервалов.
 
{\bf Доказательство.} 
Пусть $\kb{F}\in {\cal K}(C)$.
Берем $\ab =[\un{a} ,\ov{a}  ] \in\kb{IF}$,   
$\bb=[\un{b} ,\ov{b}] \in\kb{IF}$ и покажем, что $\ab + \bb $ -- интервал. Так как всегда 
$\ab + \bb \subseteq [\un{a}+\un{b} ,\ov{a}+\ov{b}]$, то достаточно показать обратное включение, то есть
$[\un{a}+\un{b} ,\ov{a}+\ov{b}]\subseteq \ab + \bb $. 
\newline
\indent
Пусть $x\in [\un{a}+\un{b} ,\ov{a}+\ov{b}]$. Тогда $0\leq x-(\un{a}+\un{b})\leq \ov{a}+\ov{b}-(\un{a}+\un{b})$ и существует 
$\tau\in {\bf\Lambda} (\kb{F})$ такое, что $\tau (\ov{a}+\ov{b}-(\un{a}+\un{b}))=x-(\un{a}+\un{b})$.
Обозначим $x_1=\un{a}+\tau (\ov{a}-\un{a})$ и $x_2=\un{b}+\tau (\ov{b}-\un{b})$. 
Тогда очевидно, что $x_1 \in [\un{a} ,\ov{a} ]$, $x_2 \in [\un{b} ,\ov{b} ]$ и $x_1 + x_2 =x$, т. e. $x\in \ab + \bb $.
Лемма доказана.

Некоторые другие характеризации пространств, обладающих свойствами $C$ и $D$, можно найти в работах \cite{lak1, lak2}. 
В частности, там доказано следующее утверждение, показывающее, что классы ${\cal K}(C)$ и ${\cal K}(D)$ достаточно обширные.

{\bf Теорема 4.} Класс ${\cal K}(C)\bigcap {\cal K}(D)$ содержит в себе
\newline
\indent
(а) все конечномерные векторные решетки;
\newline
\indent
(б) все векторные решетки, являющиеся $K_\sigma$-пространствами.

Далее будем представлять интервалы и интервальные операторы в цен\-траль\-но-сим\-мет\-рич\-ной форме. 
Если $\bb=[\un{b} ,\ov{b}] \in\kb{IF}$ -- интервал в $\kb{F}$, $\A=[\un{A} ,\ov{A}] \in\kb{IL(E,F)}$ -- интервальный оператор, 
то, полагая 
$$
\check{b}=\frac{1}{2}(\ov b + \un b),\;  \delta=\frac{1}{2}(\ov b - \un b), \;\;
\check{A}=\frac{1}{2}(\ov A + \un A),\;  \Delta=\frac{1}{2}(\ov A - \un A),
$$
получаем
$$
\bb=\left[\check{b} -\delta, \check{b} +\delta \right], \;\; 
\A=\left[\check{A} -\Delta, \check{A} +\Delta \right].
$$
При этом очевидно, что $\delta\in\kb{F}_+$ и $\Delta\in\kb{L(E,F)}_+$.
Отметим также, что в литературе имеются и другие обозначения для этих величин которые мы использовали во введении, а именно, 
${\rm mid}(\bb)=\frac{1}{2}(\ov b + \un b)$ и 
${\rm rad}(\bb)=\frac{1}{2}(\ov b - \un b)=\frac{1}{2}{\rm wid}(\bb)$ (аналогично для интервальных 
операторов \cite{Neum, sarb}).

Следующая теорема также имеется в \cite{lak1,  lak2}, однако ввиду ее важности для дальнейшего мы приведем ее  с доказательством.

{\bf Теорема 5.}  Пусть $\kb{E}\in {\cal K}(D)$, $\kb{F}\in {\cal K}(C)$, 
$\A=[\check{A}-\Delta,\check{A}+\Delta]\in \kb{IL(E,F)}$ и $x\in \kb{E}$. 
\newline
\indent
Тогда $\A x$ -- интервал в $\kb{F}$ и 
$$
\A x = [\check{A} x-\Delta |x|,\check{A} x+\Delta |x|].
\eqno (14)
$$

{\bf Доказательство.} Покажем, что $\A x \subseteq [\check{A} x-\Delta |x|,\check{A} x+\Delta |x|]$. 
Если $A\in\A$, то $A=\check{A}+\tilde{A}$ для некоторого $\tilde{A}\in [-\Delta ,\Delta]$. Тогда из неравенств
$$
Ax=\check{A}x +\tilde{A}x=\check{A}x +\tilde{A}x^+ -\tilde{A}x^-\leq \check{A}x +\Delta x^+ +\Delta x^- =
\check{A}x +\Delta (x^+ +x^-)=\check{A}x +\Delta |x|,
$$
$$
Ax=\check{A}x +\tilde{A}x=\check{A}x +\tilde{A}x^+ -\tilde{A}x^-\geq \check{A}x -\Delta x^+ -\Delta x^- =
\check{A}x -\Delta (x^+ +x^-)=\check{A}x -\Delta |x|,
$$
получаем, что $Ax\in [\check{A} x-\Delta |x|,\check{A} x+\Delta |x|]$ и, следовательно,
\newline
$\A x \subseteq [\check{A} x-\Delta |x|,\check{A} x+\Delta |x|]$.

Покажем обратное включение. Пусть $y\in [\check{A} x-\Delta |x|,\check{A} x+\Delta |x|] $. Тогда очевидно, что 
$$
0\leq y-\check{A} x+\Delta |x|\leq 2\Delta |x|.
$$

По свойству $C$ для $\kb{F}$ 
найдется мультипликатор  $\tau\in {\bf\Lambda} (\kb{F})$ такой, что
$$
y-\check{A} x+\Delta |x|= 2\tau\Delta |x|,
$$
т.е. $y=\check{A} x+(2\tau -I_{\kb{F}})\Delta |x|.$

Далее, по свойству  $D$ для $\kb{E}$ найдется мультипликатор  $\zeta\in {\bf\Lambda} (\kb{E})$ такой, что
$\zeta x^+ =x^+$ и $\zeta x^-=0$. Тогда для $\zeta$ выполняются равенства
$$
(2\zeta -I_{\kb{E}})x=2\zeta (x^+ -x^-) -(x^+ -x^-)=2\zeta x^+ -2\zeta x^--x^+ +x^-=2x^+ -x^+ +x^-=|x|.
$$

Поэтому 
$$
y=\check{A} x+(2\tau -I_{\kb{F}})\Delta |x|=\check{A} x+(2\tau -I_{\kb{F}})\Delta (2\zeta -I_{\kb{E}})x =
(\check{A} +(2\tau -I_{\kb{F}})\Delta (2\zeta -I_{\kb{E}}))x.
$$

Обозначим $\tilde{A} =(2\tau -I_{\kb{F}})\Delta (2\zeta -I_{\kb{E}})$. Тогда из следующих двух легко проверяемых тождеств 
$$
\Delta + \tilde{A} =\Delta + (2\tau -I_{\kb{F}})\Delta (2\zeta -I_{\kb{E}})=2\tau\Delta \zeta 
+2(I_{\kb{F}}-\tau)\Delta (I_{\kb{E}} - \zeta), 
$$
$$
\Delta - \tilde{A} =\Delta - (2\tau -I_{\kb{F}})\Delta (2\zeta -I_{\kb{E}})=2(I_{\kb{F}}-\tau)\Delta \zeta 
+2\tau\Delta (I_{\kb{E}} - \zeta), 
$$
и очевидных неравенств (вытекающих из того, что $0_{\kb{F}}\leq\tau\leq I_{\kb{F}}, \; 0_{\kb{E}}\leq\zeta\leq I_{\kb{E}}$ 
и $\Delta\geq 0$)
$$
2\tau\Delta \zeta +2(I_{\kb{F}}-\tau)\Delta (I_{\kb{E}} - \zeta)\geq 0,\;\;
2(I_{\kb{F}}-\tau)\Delta \zeta +2\tau\Delta (I_{\kb{E}} - \zeta)\geq 0,
$$
получаем, что
$$
\Delta + \tilde{A} \geq 0,\;\;\Delta - \tilde{A}\geq 0,
$$
т.е. $-\Delta \leq \tilde{A}\leq\Delta $. 

Следовательно, $\check{A} +\tilde{A}\in [\check{A}-\Delta,\check{A}+\Delta] $ и $y=(\check{A} +\tilde{A} )x\in \A x$. Теорема доказана.

{\bf Замечание 4.}  В случае, когда $\kb{E}=\kb{R}^n$, $\kb{F}=\kb{R}^m$ с покоординатными
упорядочениями, данная теорема повторяет, по существу, известную
теорему  Оеттли-Прагера \cite{oetpr, Oet}. Ее формулировка в виде формулы (14) имеется в \cite[Proposition 2.27, стр. 62]{FNR} 
и в \cite[Предложение 2.2.4, стр. 94-95]{sarb}.

Из этой теоремы и теоремы 3 получаем следующее бескванторное описание множества  решений уравнения (9a).
При этом считаем, что интервальные операторы и интервалы представлены в центрально-симметричной форме, т.е.

$\A_i^\prime =[\check{A}_i^\prime -\Delta_i^\prime ,\check{A}_i^\prime +\Delta_i^\prime ],\;\;
\A_i^{\prime\prime} = [\check{A}_i^{\prime\prime} -\Delta_i^{\prime\prime} ,\check{A}_i^{\prime\prime} +\Delta_i^{\prime\prime}] \in \kb{IL(E,F)},$

$\bb_i^\prime =[\check{b}_i^\prime -\delta_i^\prime ,\check{b}_i^\prime +\delta_i^\prime ] ,\;\;
\bb_i^{\prime\prime} =[\check{b}_i^{\prime\prime} -\delta_i^{\prime\prime} ,\check{b}_i^{\prime\prime} +\delta_i^{\prime\prime}] \in \kb{IF},\;\;\; i=\overline{1,\kappa}$.

{\bf Теорема 6.} Пусть $\kb{E}\in {\cal K}(D)$, $\kb{F}\in {\cal K}(C)$,   
заданы два кортежа интервальных операторов 
$\A_\forall =\left\langle \A_1^\prime , \ldots ,\A_{\kappa}^\prime \right\rangle$, 
$\A_\exists =\left\langle \A_1^{\prime\prime} , \ldots ,\A_{\kappa}^{\prime\prime} \right\rangle$, где
$\A_i^\prime ,\A_i^{\prime\prime} \in \kb{IL(E,F)}$, $i=\overline{1,\kappa}$ и  два кортежа интервалов 
$\bb_\forall =\left\langle \bb_1^\prime , \ldots ,\bb_{\kappa}^\prime \right\rangle$, 
$\bb_\exists =\left\langle \bb_1^{\prime\prime} , \ldots ,\bb_{\kappa}^{\prime\prime} \right\rangle$, где 
$\bb_i^\prime ,\bb_i^{\prime\prime} \in \kb{IF}$, $i=\overline{1,\kappa}$.
\newline
\indent
Тогда $x\in \varXi_{iq}^{\kappa}(\A_\forall , \A_\exists , \bb_\forall , \bb_\exists )$, если и только если выполняются 
следующие условия:
$$
\begin{array}{c}
-\sum\limits_{i=1}^{\kappa}\left(\Delta_i^{\prime\prime}|x| + \delta_i^{\prime\prime} \right) + 
\sum\limits_{i=1}^{\kappa}\left(\Delta_i^\prime |x| + \delta_i^\prime\right) \leq 
\sum\limits_{i=1}^{\kappa}\left((\check{A}_i^\prime + \check{A}_i^{\prime\prime}) x 
-(\check{b}_i^\prime + \check{b}_i^{\prime\prime})\right) \leq \\ \\
\leq\sum\limits_{i=1}^{\kappa}\left(\Delta_i^{\prime\prime}|x| + \delta_i^{\prime\prime} \right) - 
\sum\limits_{i=1}^{\kappa}\left(\Delta_i^\prime |x| + \delta_i^\prime\right),
\end{array} 
\eqno (15)
$$
$$
\sum\limits_{i=1}^{l}\left(\Delta_i^\prime |x| + \delta_i^\prime\right) \leq 
\sum\limits_{i=1}^{l}\left(\Delta_i^{\prime\prime}|x| + \delta_i^{\prime\prime} \right) , 
\;\;\;\; l= \ov {1, \kappa -1}.
\eqno (16)
$$

{\bf Доказательство.} Из теоремы 5 и леммы 3 следует, что выполняются условия теоремы 3. Поэтому 
$x\in \varXi_{iq}^{\kappa}(\A_\forall , \A_\exists , \bb_\forall , \bb_\exists )$ тогда и только тогда, когда выполняются условия (11), (12) теоремы 3. Кроме того, очевидно, что 
${\rm wid}(\A_i^\prime x - \bb_i^\prime)=2(\Delta_i^\prime |x| + \delta_i^\prime)$ и ${\rm wid}(\bb_i^{\prime\prime} - \A_i^{\prime\prime}  x)=2(\Delta_i^{\prime\prime}|x| + \delta_i^{\prime\prime})$ 
для всех $l= \ov {1, \kappa -1}$. Поэтому условие (12) превращается в условие (16).

Далее преобразуем условие (11), используя формулу (14) теоремы 5. Вводя обозначения 
$$
\begin{array}{c}
\check{A}^\prime =\sum\limits_{i=1}^{l}\check{A}_i^\prime ,\;\;\;
\Delta^\prime =\sum\limits_{i=1}^{l}\Delta_i^\prime ,\;\;\;
\check{A}^{\prime\prime} =\sum\limits_{i=1}^{l}\check{A}_i^{\prime\prime} ,\;\;\;
\Delta^{\prime\prime} =\sum\limits_{i=1}^{l}\Delta_i^{\prime\prime} ,\;\;\; \\
\check{b}^\prime =\sum\limits_{i=1}^{l}\check{b}_i^\prime ,\;\;\;
\delta^\prime =\sum\limits_{i=1}^{l}\delta_i^\prime ,\;\;\;
\check{b}^{\prime\prime} =\sum\limits_{i=1}^{l}\check{b}_i^{\prime\prime} ,\;\;\;
\delta^{\prime\prime} =\sum\limits_{i=1}^{l}\delta_i^{\prime\prime} ,
\end{array}
\eqno (17)
$$
и используя теорему о сумме интервалов, получаем, что условие (11) эквивалентно следующему включению:
$$
\begin{array}{c}
\left[\check{A}^\prime x -\Delta^\prime |x|-\check{b}^\prime -\delta^\prime ,
\check{A}^\prime x +\Delta^\prime |x|-\check{b}^\prime +\delta^\prime\right]\subseteq \\ \\
\subseteq
\left[\check{b}^{\prime\prime} -\delta^{\prime\prime} - \check{A}^{\prime\prime} x -\Delta^{\prime\prime} |x|,
\check{b}^{\prime\prime} +\delta^{\prime\prime} - \check{A}^{\prime\prime} x +\Delta^{\prime\prime} |x|\right],
\end{array} 
$$
которое, в свою очередь,  эквивалентно двум неравенствам
$$
\begin{array}{c}
\check{b}^{\prime\prime} -\delta^{\prime\prime} - \check{A}^{\prime\prime} x -\Delta^{\prime\prime} |x|\leq
\check{A}^\prime x -\Delta^\prime |x|-\check{b}^\prime -\delta^\prime , \;\\ \\
\check{A}^\prime x +\Delta^\prime |x|-\check{b}^\prime +\delta^\prime\leq
\check{b}^{\prime\prime} +\delta^{\prime\prime} - \check{A}^{\prime\prime} x +\Delta^{\prime\prime} |x|.
\end{array}
$$

Эти неравенства с помощью очевидных преобразований приводятся к следующему виду: 
$$
-(\Delta^{\prime\prime} |x|+\delta^{\prime\prime})+(\Delta^\prime |x|+\delta^\prime) \leq
(\check{A}^\prime +\check{A}^{\prime\prime})x-(\check{b}^\prime +\check{b}^{\prime\prime}) \leq
(\Delta^{\prime\prime} |x|+\delta^{\prime\prime})-(\Delta^\prime |x|+\delta^\prime) .
$$

Подставляя в последнее неравенство обозначения (17), получаем условие (15). Теорема доказана.

В случае, когда $\kb{F}$ является векторной решеткой, условие (15) можно несколько упростить.

{\bf Теорема 7.} Пусть выполняются условия теоремы 6 и, кроме того, $\kb{F}$ является векторной решеткой.
Тогда $x\in \varXi_{iq}^{\kappa}(\A_\forall , \A_\exists , \bb_\forall , \bb_\exists )$, если и только если выполняется 
условие (16) и  следующее условие:
$$
\left|\sum\limits_{i=1}^{\kappa}\left((\check{A}_i^\prime + \check{A}_i^{\prime\prime}) x 
-(\check{b}_i^\prime + \check{b}_i^{\prime\prime})\right)\right| +
\sum\limits_{i=1}^{\kappa}\left(\Delta_i^\prime |x| + \delta_i^\prime\right)\leq 
\sum\limits_{i=1}^{\kappa}\left(\Delta_i^{\prime\prime}|x| + \delta_i^{\prime\prime} \right).  
\eqno (15^\prime)
$$

Для доказательства достаточно заметить, что если $\kb{F}$ является векторной решеткой, то  в ней для любых элементов 
$y,z\in\kb{F}$ неравенство $-z\leq y\leq z$ эквивалентно неравенству $|y|\leq z$.

{\bf Замечание 5.}  Отметим, что формулы теоремы 7 можно получить и с помощью детального анализа процедуры 
элиминации кван\-то\-ров из работы В. Крейновича \cite{Krei} (в конечномерном случае $m$ уравнений от $n$ переменных и с покоординатным упорядочением), примененной им для построения алгоритма порядка $O(m\cdot n)$ проверки принадлежности вектора множеству решений  интервально-кванторной линейной системы уравнений. Обратно, с помощью  формул теоремы 7, такой алгоритм строится очевидным образом.

Если  $\kb{E}=\R^n$ и $\kb{F}=\R^m$ с покоординатными упорядочиваниями, то из теоремы 7 и предложения 1 получаем, что 
класс множеств решений уравнения (9) также совпадает с 
классом AE-решений. Более точно верно следующее утверждение.

{\bf Предложение 2.} Пусть $\kb{E}=\R^n$ и $\kb{F}=\R^m$ с покоординатными упорядочиваниями. Тогда 
для любых кортежей интервальных $m\times n$-матриц 
\linebreak
$\A_\forall =\left\langle \A_1^\prime , \ldots ,\A_{\kappa}^\prime \right\rangle$,
$\A_\exists =\left\langle \A_1^{\prime\prime} , \ldots ,\A_{\kappa}^{\prime\prime} \right\rangle$
и  кортежей интервальных $m$-векторов 
\linebreak
$\bb_\forall =\left\langle \bb_1^\prime , \ldots ,\bb_{\kappa}^\prime \right\rangle$, 
$\bb_\exists =\left\langle \bb_1^{\prime\prime} , \ldots ,\bb_{\kappa}^{\prime\prime} \right\rangle$
существуют интервальная 
$(\kappa m)\times n$-матрица $\tilde{\A}$, интервальный $(\kappa m)$-вектор $\tilde{\bb}$ и кванторные 
$(\kappa m)\times n$-матрица $\alpha$ и $(\kappa m)$-вектор $\beta$ такие, что выполняется равенство
$$\varXi_{\alpha\beta}(\tilde{\A}, \tilde{\bb}) = 
\varXi_{iq}^{\kappa}(\A_\forall , \A_\exists , \bb_\forall , \bb_\exists ).
$$

{\bf Доказательство.} Представим неравенства (16), $(15^\prime)$ в виде одной системы неравенств с модулями $(lm)$, полагая
$$
C= 
\left( 
\begin{array}{c}
C_1 \\
\vdots \\
C_{\kappa}
\end{array}
\right) ,
D= 
\left( 
\begin{array}{c}
D_1 \\
\vdots \\
D_{\kappa}
\end{array}
\right) ,
c= 
\left( 
\begin{array}{c}
c_1 \\
\vdots \\
c_{\kappa}
\end{array}
\right) ,
d= 
\left( 
\begin{array}{c}
d_1 \\
\vdots \\
d_{\kappa}
\end{array}
\right) ,
$$
где $C_l =0$, $c_l =0$ при $l=\overline{1,\kappa -1}$, 
$C_{\kappa} = \sum\limits_{s=1}^{\kappa}\left(\check{A}_s^\prime + \check{A}_s^{\prime\prime}\right)$, 
$c_{\kappa} = \sum\limits_{s=1}^{\kappa}\left(\check{b}_s^\prime + \check{b}_s^{\prime\prime}\right)$ 
и  
$D_l =$ $= \sum\limits_{s=1}^{l}\left(\Delta_s^{\prime\prime} -\Delta_s^\prime \right)$, 
$d_l = \sum\limits_{s=1}^{l}\left(\delta_s^{\prime\prime} -\delta_s^\prime \right)$ для всех $l=\overline{1,\kappa}$.

Тогда, определяя матрицы $\tilde{\A}$, $\alpha$ и векторы $\tilde{\bb}$, $\beta$ по этим $C,\; c, \; D, \; d$ 
так же, как в предложении 1, получаем, что  
$\varXi_{\alpha\beta}(\tilde{\A}, \tilde{\bb}) = 
\varXi_{iq}^{\kappa}(\A_\forall , \A_\exists , \bb_\forall , \bb_\exists ).
$
Предложение 2 доказано.

\begin{center}
{\bf Заключение }
\end{center}

В работе рассматриваются множества решений  интервально-кванторных линейных систем уравнений и их обобщений.
С помощью выделения в кванторной приставке $\forall\exists$-блоков получен  способ приведения интервально-кванторных линейных систем уравнений к некоторому каноническому виду, в котором кванторы всеобщности и существования строго чередуются. Это также позволило несколько обобщить  понятие интервально-кванторных линейных систем уравнений и перенести его на произвольные упорядоченные векторные пространства. Получено бескванторное описание этих множеств  как в интервальной арифметике -- теорема 3 (что обобщает соответствующие результаты H. Beeck'а \cite{Beeck} и С.П. Шарого \cite{sar4}), так и 
в виде разрешимости систем  линейных неравенств с модулями -- теоремы 6 и 7 (что обобщает  результат Оеттли—Прагера для объединенного множества решений \cite{oetpr} и результат И. Рона для AE-решений \cite{R1, R2}). 
Используя бескванторное описание в виде разрешимости систем  линейных неравенств с модулями, показано, что 
в конечномерном случае с  покоординатным упорядочиванием 
класс множеств решений обобщенных интервально-кванторных уравнений совпадает с классом AE-решений. 
Это позволяет применять уже созданные методы (в том числе и численные) для их исследования.


\end{document}